\documentclass[11pt]{amsart}

\usepackage{graphicx}
\usepackage{bbm}
\usepackage{hyperref}
\usepackage{breakurl}
\usepackage{amsfonts,amsmath,amssymb}
\usepackage{xcolor}
\usepackage{tikz}
\usetikzlibrary{3d}
\usetikzlibrary{calc}
\usepackage[justification=raggedright,font=footnotesize]{caption}
\usepackage[labelfont=rm,font=footnotesize]{subcaption}
\usepackage{booktabs,tabularx,colortbl}
\usepackage{dcolumn}
\newcolumntype{d}[1]{D{,}{.}{-1}}

\usepackage{listings}
\newcommand{\R}{\mathbbm R}
\newcommand{\Z}{\mathbbm Z}
\newcommand{\conv}{\mathsf{conv}}
\newcommand{\1}{\boldsymbol{1}}
\newcommand{\polydb}{\texttt{polyDB}}
\newcommand{\polymake}{\texttt{polymake}}
\DeclareMathOperator{\sbip}{SBipyr}
\definecolor{lightgray}{rgb}{.8,.8,.8}
\newtheorem{theorem}{Theorem}[section]
\newtheorem{proposition}[theorem]{Proposition}

\keywords{Database, polyDB, Discrete Geometry, Lattice Polytopes, Reflexive Polytopes, Smooth Polytopes}

\subjclass{ 52-04, 52B20}

\title{\polydb: A Database for Polytopes and related Objects}
\author{Andreas Paffenholz}
\address{TU Darmstadt, Dolivostr.\ 15, 64293 Darmstadt, Germany}
\thanks{The author was partially supported by the Priority Program 1489 of the German Research Council (DFG)}
\email{paffenholz@opt.tu-darmstadt.de}

\begin{document}

\begin{abstract}
  \polydb\ is a database for discrete geometric
  objects independent of a particular software. The database is
  accessible via web and an interface from the software package
  \polymake. It contains various datasets from the area of lattice
  polytopes, combinatorial polytopes, matroids and tropical geometry.
\newline\indent
  In this short note we introduce the structure of the database and
  explain its use with a computation of the free sums and certain skew
  bipyramids among the class of smooth Fano polytopes in dimension up
  to $8$.
\end{abstract}

\newsavebox{\dbentry}
\begin{lrbox}{\dbentry}
  \begin{lstlisting}
{
    "_id" : "F.2D.3",
    "DIM":2,
    "FACETS" : [[1,0,1],[1,0,-1],[1,1,0],[1,-1,-1],[1,-1,0]],
    "VERTICES":[[1,-1,-1],[1,-1,1],[1,0,1],[1,1,0],[1,1,-1]],
    "F_VECTOR" : [5,5],
    "EHRHART_POLYNOMIAL_COEFF":["1","7/2","7/2"],
    "H_STAR_VECTOR":[1,5,1],
    "CENTROID":["1","-2/21","-2/21"],
    "N_LATTICE_POINTS":8,
    "NORMAL" : "true",
    "VERY_AMPLE" : "true",
    "LATTICE_VOLUME":7,
    "polyDB" : { 
                 [...]
               }
}  
\end{lstlisting}
\end{lrbox}

\newsavebox{\productscript}
\begin{lrbox}{\productscript}
  \begin{lstlisting}
use application "polytope";

sub identify_smooth_fano_in_polydb {
  my $p = shift;
  my $d = $p->DIM;
  my $nlp = new Int($p->N_LATTICE_POINTS);
  my $parray=db_query({"DIM"=>$d, "N_VERTICES"=>$p->N_VERTICES, "N_FACETS"=>$p->N_FACETS, 
    "N_LATTICE_POINTS"=>$nlp, }, db=>"LatticePolytopes", collection=>"SmoothReflexive");
  foreach my $c ( @$parray ) {
    if ( lattice_isomorphic_smooth_polytopes($c,$p) ) { return $c->name; }
  }
  die "polytope not found\n";
}

sub all_free_sums_in_dim { 
  my ($d,%options) = @_;
  my $list;
  if ( $options{"splitinfo"} ) {
    $list = new Map<String,Set<Pair<String,String> > >;
  } else { 
    $list = new Set<String>; 
  }
  my $cur_options = { db=>"LatticePolytopes", collection=>"SmoothReflexive" };
  foreach my $n (1..$d/2) {
    my $cur1=db_cursor({"DIM"=>$n}, $cur_options);
    while ( !$cur1->at_end() ) {
      my $c1 = $cur1->next();
      my $cur2=db_cursor({"DIM"=>$d-$n}, $cur_options );
      while ( !$cur2->at_end() ) {
        my $c2 = $cur2->next();
        my $name = identify_smooth_fano_in_polydb(product($c1,$c2));
        if ( $options{"splitinfo"} ) {
          my $split = new Pair<String,String>($c1->name,$c2->name);
          $list->{$name} += $split;
        } else {
          $list += $name;
  } } } }
  return $list;
}
\end{lstlisting}
\end{lrbox}

\newcommand{\smoothFano}{
    \begin{subfigure}{.18\linewidth}
      \begin{tikzpicture}[scale=0.9]
        \tikzstyle{edge} = [draw,thick,-,black]
		  
        \foreach \x in {-1,0,1} \foreach \y in {-1,0,1} \fill[gray]
        (\x,\y) circle (1.5pt); ;
		
        \coordinate (v0) at (0,0); \coordinate (e1) at (1,0);
        \coordinate (e2) at (0,1); \coordinate (-e1+e2) at (-1,1);
        \coordinate (-e2+e1) at (1,-1); \coordinate (-e1) at (-1,0);
        \coordinate (-e2) at (0,-1);
		  
        \fill[opacity=.2] (e1) -- (e2) -- (-e1+e2) -- (-e1) -- (-e2) --
        (-e2+e1) -- (e1);
        \draw[edge] (e1) -- (e2) -- (-e1+e2) -- (-e1) -- (-e2) --
        (-e2+e1) -- (e1);
		  
        \foreach \point in {e1,e2,-e1+e2,-e2+e1,-e1,-e2} \fill[black]
        (\point) circle (2pt);
		
      \end{tikzpicture}
      \caption{$P_6$}
    \end{subfigure}
    \begin{subfigure}{.18\linewidth}
        	\begin{tikzpicture}[scale=0.9]
        	  \tikzstyle{edge} = [draw,thick,-,black]
		
        	  \foreach \x in {-1,0,1}
        	    \foreach \y in {-1,0,1}
        	       \fill[gray] (\x,\y) circle (1.5pt); ;
		
        	  \coordinate (v0) at (0,0);
        	  \coordinate (e1) at (1,0);
        	  \coordinate (e2) at (0,1);
        	  \coordinate (-e1+e2) at (-1,1);
        	  \coordinate (-e2+e1) at (1,-1);
        	  \coordinate (-e2) at (0,-1);
		  
        	  \fill[opacity=.2] (e1) -- (e2) -- (-e1+e2) -- (-e2) -- (-e2+e1) -- (e1);
        	  \draw[edge] (e1) -- (e2) -- (-e1+e2) -- (-e2) -- (-e2+e1) -- (e1);
		  
        	  \foreach \point in {e1,e2,-e1+e2,-e2+e1,-e2}
        	    \fill[black] (\point) circle (2pt);
		  
        	\end{tikzpicture}
      \caption{$P_5$}
    \end{subfigure}
    \begin{subfigure}{.18\linewidth}
        	\begin{tikzpicture}[scale=0.9]
        	  \tikzstyle{edge} = [draw,thick,-,black]
		
        	  \foreach \x in {-1,0,1}
        	    \foreach \y in {-1,0,1}
        	       \fill[gray] (\x,\y) circle (1.5pt); ;
		
        	  \coordinate (v0) at (0,0);
        	  \coordinate (e1) at (1,0);
        	  \coordinate (e2) at (0,1);
        	  \coordinate (-e1) at (-1,0);
        	  \coordinate (-e2) at (0,-1);
		
        	  \fill[opacity=.2] (e1) -- (e2) -- (-e1) -- (-e2) -- (e1);
        	  \draw[edge] (e1) -- (e2) -- (-e1) -- (-e2) -- (e1);
		  
        	  \foreach \point in {e1,e2,-e1,-e2}
        	    \fill[black] (\point) circle (2pt);
        	\end{tikzpicture}
      \caption{$P_{4a}$}
    \end{subfigure}
    \begin{subfigure}{.18\linewidth}
        	\begin{tikzpicture}[scale=0.9]
        	  \tikzstyle{edge} = [draw,thick,-,black]
		
        	   \foreach \x in {-1,0,1}
        	    \foreach \y in {-1,0,1}
        	       \fill[gray] (\x,\y) circle (1.5pt); ;
		
        	  \coordinate (v0) at (0,0);
        	  \coordinate (e1) at (1,0);
        	  \coordinate (e2) at (0,1);
        	  \coordinate (-e2+e1) at (1,-1);
        	  \coordinate (-e1) at (-1,0);
		  
        	  \fill[opacity=.2] (e1) -- (e2) -- (-e1) -- (-e2+e1) -- (e1);
        	  \draw[edge] (e1) -- (e2) -- (-e1) -- (-e2+e1) -- (e1);
		  
        	  \foreach \point in {e1,e2,-e2+e1,-e1}
        	    \fill[black] (\point) circle (2pt);
        	\end{tikzpicture}
      \caption{$P_{4b}$}
    \end{subfigure}
    \begin{subfigure}{.18\linewidth}
        	\begin{tikzpicture}[scale=0.9]
        	  \tikzstyle{edge} = [draw,thick,-,black]
		  
        	  \foreach \x in {-1,0,1}
        	    \foreach \y in {-1,0,1}
        	       \fill[gray] (\x,\y) circle (1.5pt); ;
		
        	  \coordinate (v0) at (0,0);
        	  \coordinate (e1) at (1,0);
        	  \coordinate (e2) at (0,1);
        	  \coordinate (-e1-e2) at (-1,-1);
		
        	  \fill[opacity=.2] (e1) -- (e2) -- (-e1-e2) -- (e1);
        	  \draw[edge] (e1) -- (e2) -- (-e1-e2) -- (e1);
		  
        	  \foreach \point in {e1,e2,-e1-e2}
        	    \fill[black] (\point) circle (2pt);
        	\end{tikzpicture}
      \caption{$P_3$}
    \end{subfigure}
}

\newcommand{\properBip}{
		\begin{tikzpicture}[x  = {(0cm,1.2cm)},
		                    y  = {(2cm,0cm)},
		                    z  = {(.9cm,.4cm)},
		                    scale = .95,
		                    color = {lightgray}]
		\tikzset{facestyle/.style={fill=black!20, opacity=.4,draw=black,line join=round}}
		\tikzset{factorstyle/.style={fill=black!70, opacity=.8,draw=black,line join=round}}
		\tikzset{factoredgestyle/.style={draw=black,thick}}
		
		\tikzstyle{every label}=[black]
		
		\draw[white] (1.5,0,0) -- (-1.5,0,0);
		
		  \draw[facestyle] (1,0,0) -- (0,1,0) -- (0,0,1) -- (1,0,0) -- cycle ;
		  \draw[facestyle] (1,0,0) -- (0,-1,1) -- (0,0,1) -- (1,0,0) -- cycle ;
		  \draw[facestyle] (1,0,0) -- (0,-1,1) -- (0,-1,0) -- (1,0,0) -- cycle ;
		
		  \draw[facestyle] (-1,0,0) -- (0,1,0) -- (0,0,1) -- (-1,0,0) -- cycle ;
		  \draw[facestyle] (-1,0,0) -- (0,-1,1) -- (0,0,1) -- (-1,0,0) -- cycle ;
		  \draw[facestyle] (-1,0,0) -- (0,-1,1) -- (0,-1,0) -- (-1,0,0) -- cycle ;
		
		    \draw (0,0,1) node [label=right:\raisebox{3ex}{$\scriptstyle e_3$}] {};
		    \fill[black] (0,0,1) circle (2pt);
		    \fill[black] (0,-1,1) circle (2pt);

                  \draw[factoredgestyle] (-1,0,0) -- (0,0,0);
                  \draw[factorstyle] (0,-1,0) -- (0,0,-1) -- (0,1,-1) -- (0,1,0) -- (0,0,1) -- (0,-1,1) -- cycle ;
		  \fill[black] (0,0,0) circle (2pt);
                  \draw[factoredgestyle] (1,0,0) -- (0,0,0);

		  \draw[facestyle] (1,0,0) -- (0,-1,0) -- (0,0,-1) -- (1,0,0) -- cycle ;
		  \draw[facestyle] (1,0,0) -- (0,1,-1) -- (0,0,-1) -- (1,0,0) -- cycle ;
		  \draw[facestyle] (1,0,0) -- (0,1,-1) -- (0,1,0) -- (1,0,0) -- cycle ;

		  \draw[facestyle] (-1,0,0) -- (0,-1,0) -- (0,0,-1) -- (-1,0,0) -- cycle ;
		  \draw[facestyle] (-1,0,0) -- (0,1,-1) -- (0,0,-1) -- (-1,0,0) -- cycle ;
		  \draw[facestyle] (-1,0,0) -- (0,1,-1) -- (0,1,0) -- (-1,0,0) -- cycle ;
		
		  \draw (1,0,0) node [label=above:$\scriptstyle e_1$] {};
		  \draw (-1,0,0) node [label=below:$\scriptstyle -e_1$] {};
		  \draw (0,1,0) node [label=right:$\scriptstyle e_2$] {};
		  \fill[black] (1,0,0) circle (2pt);
		  \fill[black] (-1,0,0) circle (2pt);
		  \fill[black] (0,1,0) circle (2pt);
		  \fill[black] (0,-1,0) circle (2pt);
		  \fill[black] (0,1,-1) circle (2pt);
		  \fill[black] (0,0,-1) circle (2pt);
		
		\end{tikzpicture} 
}

\newcommand{\skewBip}{
		\begin{tikzpicture}[x  = {(0cm,1.2cm)},
		                    y  = {(2cm,0cm)},
		                    z  = {(.9cm,.4cm)},
		                    scale = .95,
		                    color = {lightgray}]
		\tikzset{facestyle/.style={fill=black!20, opacity=.4,draw=black,line join=round}}
		\tikzset{factorstyle/.style={fill=black!70, opacity=.8,draw=black,line join=round}}
		\tikzset{factoredgestyle/.style={draw=black,thick}}
		\tikzstyle{vertex}=[circle,minimum size=3pt,inner sep=0pt, fill=black]
		\tikzstyle{every label}=[black]
		
		\draw[white] (1.5,0,0) -- (-1.5,0,0);
		
		  \draw[facestyle] (-1,0,0) -- (0,1,0) -- (0,0,1) -- (-1,0,0) -- cycle ;
		  \draw[facestyle] (-1,0,0) -- (0,-1,1) -- (0,0,1) -- (-1,0,0) -- cycle ;
		  \draw[facestyle] (-1,0,0) -- (0,-1,1) -- (0,-1,0) -- (-1,0,0) -- cycle ;
		
		  \draw[facestyle] (1,0,1) -- (0,1,0) -- (0,0,1) -- (1,0,1) -- cycle ;
		  \draw[facestyle] (1,0,1) -- (0,-1,1) -- (0,0,1) -- (1,0,1) -- cycle ;
		  \draw[facestyle] (1,0,1) -- (0,-1,1) -- (0,-1,0) -- (1,0,1) -- cycle ;
		
		  \draw[facestyle] (1,0,1) -- (0,-1,0) -- (0,0,-1) -- (1,0,1) -- cycle ;
		
		    \draw (0,0,1) node[label=right:\raisebox{3ex}{$\scriptstyle e_3$}] {};
		    \fill[black] (0,0,1) circle (2pt);
		    \fill[black] (0,-1,1) circle (2pt);

                  \draw[factoredgestyle] (-1,0,0) -- (0,0,1/2);
                  \draw[factorstyle] (0,-1,0) -- (0,0,-1) -- (0,1,-1) -- (0,1,0) -- (0,0,1) -- (0,-1,1) -- cycle ;
		  \fill[black] (0,0,1/2) circle (2pt);
                  \draw[factoredgestyle] (1,0,1) -- (0,0,1/2);

		  \draw[facestyle] (-1,0,0) -- (0,-1,0) -- (0,0,-1) -- (-1,0,0) -- cycle ;
		  \draw[facestyle] (-1,0,0) -- (0,1,-1) -- (0,0,-1) -- (-1,0,0) -- cycle ;
		  \draw[facestyle] (-1,0,0) -- (0,1,-1) -- (0,1,0) -- (-1,0,0) -- cycle ;

		  \draw[facestyle] (1,0,1) -- (0,1,-1) -- (0,0,-1) -- (1,0,1) -- cycle ;
		  \draw[facestyle] (1,0,1) -- (0,1,-1) -- (0,1,0) -- (1,0,1) -- cycle ;
		
		  \draw (-1,0,0) node[label=below:$\scriptstyle -e_1$] {};
		  \draw (1,0,1) node[label=right:$\scriptstyle e_1+e_3$] {};
		  \draw (0,1,0) node[label=right:$\scriptstyle e_2$] {};
		  \fill[black] (-1,0,0) circle (2pt);
		  \fill[black] (1,0,1) circle (2pt);
		  \fill[black] (0,1,0) circle (2pt);
		  \fill[black] (0,-1,0) circle (2pt);
		  \fill[black] (0,1,-1) circle (2pt);
		  \fill[black] (0,0,-1) circle (2pt);
		
		\end{tikzpicture} 
}

\newcommand{\skewSimplexSum}{
		\begin{tikzpicture}[x  = {(0cm,1.2cm)},
		                    y  = {(2cm,0cm)},
		                    z  = {(.9cm,.4cm)},
		                    scale = .95,
		                    color = {lightgray}]
		\tikzset{facestyle/.style={fill=black!20, opacity=.4,draw=black,line join=round}}
		\tikzset{factorstyle/.style={fill=black!70, opacity=.8,draw=black,line join=round}}
		\tikzset{factoredgestyle/.style={draw=black,thick}}
		\tikzstyle{vertex}=[circle,minimum size=3pt,inner sep=0pt, fill=black]
		\tikzstyle{every label}=[black]
		
		\draw[white] (1.8,0,0) -- (-1.5,0,0);
		
		  \draw[facestyle] (1,0,0) -- (0,1,1) -- (0,-1,0) -- cycle ;
		  \draw[facestyle] (-1,0,0) -- (0,1,1) -- (0,-1,0) -- cycle ;

                  \draw[factoredgestyle] (-1,0,0) -- (0,0,0);
                  \draw[factorstyle] (0,-1,0) -- (0,0,-1) -- (0,1,1) -- cycle ;
		  \fill[black] (0,0,0) circle (2pt);
                  \draw[factoredgestyle] (1,0,0) -- (0,0,0);
		
		  \draw[facestyle] (0,-1,0) -- (0,0,-1) -- (-1,0,0) -- cycle ;
		  \draw[facestyle] (0,-1,0) -- (0,0,-1) -- (1,0,0) -- cycle ;
		  \draw[facestyle] (0,1,1) -- (0,0,-1) -- (-1,0,0) -- cycle ;
		  \draw[facestyle] (0,1,1) -- (0,0,-1) -- (1,0,0) -- cycle ;

		  \draw (1,0,0) node[label=above:$\scriptstyle e_1$] {};
		  \draw (0,1,1) node[label=above:$\scriptstyle e_2+e_3$] {};
		  \draw (0,-1,0) node[label=left:$\scriptstyle -e_2$] {};
                  \draw (0,0,-1) node[label=left:\raisebox{-4ex}{$\scriptstyle -e_3$}] {};
		  \fill[black] (0,0,-1) circle (2pt);
		  \fill[black] (1,0,0) circle (2pt);
		  \fill[black] (0,1,1) circle (2pt);
		  \fill[black] (0,-1,0) circle (2pt);
		  \fill[black] (-1,0,0) circle (2pt);
		
		\end{tikzpicture} 
}

\newcommand{\skewSimplexSumTwo}{
		\begin{tikzpicture}[x  = {(0cm,1.2cm)},
		                    y  = {(2cm,0cm)},
		                    z  = {(.9cm,.4cm)},
		                    scale = .95,
		                    color = {lightgray}]
		\tikzset{facestyle/.style={fill=black!20, opacity=.4,draw=black,line join=round}}
		\tikzset{factorstyle/.style={fill=black!70, opacity=.8,draw=black,line join=round}}
		\tikzset{factoredgestyle/.style={draw=black,thick}}
		\tikzstyle{vertex}=[circle,minimum size=3pt,inner sep=0pt, fill=black]
		\tikzstyle{every label}=[black]
		
		\draw[white] (1.8,0,0) -- (-1.5,0,0);

		  \draw[facestyle] (1,0,0) -- (1,1,1) -- (0,-1,0) -- cycle ;
		  \draw[facestyle] (-1,0,0) -- (1,1,1) -- (0,-1,0) -- cycle ;

                  \draw[factoredgestyle] (-1,0,0) -- (1/3,0,0);
                  \draw[factorstyle] (0,-1,0) -- (0,0,-1) -- (1,1,1) -- cycle ;
		  \fill[black] (1/3,0,0) circle (2pt);
                  \draw[factoredgestyle] (1,0,0) -- (1/3,0,0);
		
		  \draw[facestyle] (0,-1,0) -- (0,0,-1) -- (-1,0,0) -- cycle ;
		  \draw[facestyle] (0,-1,0) -- (0,0,-1) -- (1,0,0) -- cycle ;
		  \draw[facestyle] (1,1,1) -- (0,0,-1) -- (-1,0,0) -- cycle ;
		  \draw[facestyle] (1,1,1) -- (0,0,-1) -- (1,0,0) -- cycle ;

		  \draw (1,0,0) node[label=above:$\scriptstyle e_1$] {};
		  \draw (0,-1,0) node[label=left:$\scriptstyle -e_2$] {};
		  \draw (1,1,1) node[label=above:\hspace*{-.3cm}$\scriptstyle e_1+e_2+e_3$] {};
                  \draw (0,0,-1) node[label=left:\raisebox{-4ex}{$\scriptstyle -e_3$}] {};
		  \fill[black] (0,0,-1) circle (2pt);
		  \fill[black] (1,0,0) circle (2pt);
		  \fill[black] (1,1,1) circle (2pt);
		  \fill[black] (0,-1,0) circle (2pt);
		  \fill[black] (-1,0,0) circle (2pt);
		
		\end{tikzpicture} 
}

\maketitle

\section{Introduction}

In recent years availability of computational classifications of
mathematical objects has proven to be an important and valuable tool
to obtain new results, to check new ideas and to experiment with the
objects to obtain insight into their structure and directions for
further research.

We know the full list of smooth Fano polytopes (up to lattice
equivalence) up to dimension $9$ by an algorithm of \O bro
\cite{OebroPhD}, whose availability within the software package
\polymake\ has been the foundation \textit{e.g.} for counter-examples
to a conjecture of Batyrev and Selivanova~\cite{MR2842630} or the
classification of simplicial, terminal, and reflexive polytopes with
many vertices by Assarf et al.~\cite{AJP}. Availability of the same
data in Magma~\cite{magma_sw} lead to the study of the poset of
blowups by Higashitani~\cite{1503.06434} or the study of reflexive
polytopes of higher index by Kasprzyk and Nill~\cite{1107.4945}. The
classification of $0/1$-polytopes up to dimension $6$ by
Aichholzer~\cite{Aichholzer} was used in the study of permutation
polytopes by Baumeister et al.~\cite{BHNP07-1}.

We also know classifications of small oriented matroids by Miyata et
al.~\cite{MR3017917}, polytropes (Kulas, Joswig,
Tran~\cite{MR2629819,MR3201242}), and reflexive polytopes up to
dimension $4$ by Kreuzer and Skarke~\cite{cydata,1411.1418}. The
symbolic data project by Gr\"abe et al.~\cite{symbolicdata} aims to
collect data from computer algebra and make it accessible in a
structured and searchable form on their web page. The library
\textit{MIPLIB} by Koch et al.~\cite{KochEtAl2011} collects discrete
optimization problems for benchmarking of algorithms.

Most of these collections, however, cannot easily be used in a
software package. Sometimes the data is only available in text format
or, if searchable via a database, is connected to a specific software
package or lacks a proper interface at all. For example, the small
oriented matroids \cite{MR3017917}, polytropes
\cite{MR3201242,1012.3053}, or $0/1$-polytopes \cite{Aichholzer} are
available as text files, while access to the small groups library
\cite{sgl_data} is linked to GAP \cite{GAP4}.  Altman et
al.~\cite{1411.1418} have created a database for the reflexive
polytopes up to dimension $4$ computed by Kreuzer and
Skarke~\cite{cydata}, but it is currently not accessible at the link
given in the paper.

On the other hand, the Graded Rings Database~\cite{grdb} project has a
more general approach and provides data in a format both searchable
via a web interface and accessible via a programmatic interface that
can be used in software packages. It currently has a focus on data
from combinatorial commutative algebra and toric geometry.

The new database \polydb\ aims to provide searchable data from a wide
range of areas at a permanent location in an application independent
format. It allows download in text format and access from any software
package that provides an interface to the data. It is also searchable
via a web interface at
\href{http://db.polymake.org}{db.polymake.org}. Currently, one
interface to a software package is implemented, in the software
package \polymake~\cite{1408.4653,Joswig2009}. The current collection
of data is thus still inspired by the range of applications of
\polymake\ with data from combinatorial geometry, matroid theory, toric
geometry and combinatorial topology.

In the following two sections we explain the concept of the database
and introduce the interface implemented in \polymake\ to access the
data. The last section shows one application of the database and the
interface. We will show that in dimensions up to $8$ more than $80\%$
of the smooth Fano polytopes arise from lower dimensional ones as a
free sum of two lower dimensional smooth Fano polytopes or a certain
skew sum construction of a smooth Fano polytope and a simplex. We give
the count of polytopes decomposable in this way in
Table~\ref{tab:sr_prod}. With a simple extension of the scripts one
can also obtain the list of possible decompositions for each polytope.

\section{polyDB}

In this section we briefly introduce the structure of the database
\polydb\ and the data sets already contained in it.

The database \polydb\ for discrete geometric objects is based on the
open source NoSQL database MongoDB~\cite{mongodb}. It has been set up
at \href{https://db.polymake.org}{db.polymake.org}. The database
stores its data as plain JSON documents grouped into
\textit{collections} and \textit{databases} (To avoid confusion with
this and the abstract database \polydb\ we will refer to this
technical term introduced by MongoDB as a \textit{collection
  group}). We use this to group collections from the same area of
discrete geometry into a common collection group. \textit{E.g.}, the
collection group \textit{Objects in Tropical Geometry} currently
contains two collections of such objects, the small \textit{tropical
  oriented matroids} classified by Horn~\cite{hornPhD} and the
\textit{polytropes} classified by Kulas~\cite{1012.3053} and
Tran~\cite{MR3201242}. \polydb\ stores data in a plain JSON format
independent of any particular software package. See
Figure~\ref{fig:entry} for an example of an entry in the collection of
smooth reflexive polytopes.
\begin{figure}[t]
    \usebox{\dbentry} 
    \caption{An entry in the collection of smooth Fano
      polytopes. Naming of the fields is in this example taken from
      standard properties of objects in \polymake. However, there are
      no restrictions on field names.}
\label{fig:entry}
\end{figure}

Each document contains one special entry \textit{polyDB} (besides its
\textit{\_id}, which is required by MongoDB). Apart form this all
other entries and their tags can be chosen freely depending on the
data. The entry \textit{polyDB} may specify format restrictions for
the data and import or export specifications for various software
packages, separated by subfields naming the software. This section may
contain, \textit{e.g.}, information on the required version, authors
of the data, and the method to load the data into the particular
software package. Each collection group also has a separate
collection \textit{type\_information} that specifies the format of an
entry in a collection and allows to store information applicable to
all data sets in this collection, \emph{e.g.}, methods for import and
export of the data.  The web interface at
\href{http://db.polymake.org}{db.polymake.org} allows independent and
searchable access to all data sets in \polydb.

There are currently five collections, grouped into four
collection groups contained in \polydb. We give a brief introduction
to each of the collections.
\begin{itemize}
\item The collection group \textit{Lattice Polytopes} has the collection
  \textit{Smooth Reflexive Polytopes} that contains low dimensional
  smooth reflexive polytopes based on the algorithm of \O bro
  \cite{OebroPhD}. \O bro used his algorithm to compute the data up to
  dimension $8$. Later, dimension $9$ was computed with an improved
  implementation of the algorithm by Lorenz and the author. There are
  $9,060,505$ such polytopes.
\item The collection group \textit{Objects in Tropical Geometry} has two
  collections.  The collection \textit{Tropical Oriented Matroids}
  contains a list of $71$ known non-realizable tropical oriented
  matroids. This data was provided by Horn~\cite{hornPhD}.  The
  collection \textit{Full-dimensional Polytropes in TP3} contains all
  $1013$ polytropes in $3$-dimensional tropical projective space. The
  collection was generated by Constantin Fischer from data of Joswig
  and Kulas~\cite{MR2629819} and
  Tran~\cite{MR3201242}. See~\cite{webPolytropes} for a description.
\item The collection group \textit{Special Polytopes} has the collection
  \textit{Faces of Birkhoff Polytopes} which contains all $5371$
  combinatorial types of faces up to dimension $8$ of the Birkhoff
  polytope in any dimension~\cite{MR3336581}.
\item The collection group \textit{Matroids} has the collection
  \textit{Matroids on at most $12$ elements}. This collection contains
  all $32,401,446$ small matroids as computed by Miyata et
  al.~\cite{MR3017917,om_data,MR2038485}.
\end{itemize}
Further collections are in preparation.

\section{The \polymake\ interface to \polydb}

The initiative for \polydb\ was started in 2013 by Silke Horn and the
author as an \emph{extension} for the software package
\polymake~\cite{githubpolydb} with associated database. With the
latest version 3.1 of \polymake~\cite{polymakesnapshot}, released in
March 2017, the interface to the database has been turned into a
\emph{bundled extension} for \polymake\ that is directly delivered
with the software and the database has been set up as an independent
project.

However, the software package \polymake\ currently provides the only
interface for import of data into the database and methods to access
and use it for computations. Given a search query, \textit{i.e.} a
list of restrictions on the properties of an object, MongoDB allows
the retrieval of a single object satisfying the query, an array with
all objects satisfying the query or a cursor that returns objects from
the result set one after another. All three methods are also
implemented in \polymake. The implementation is based on the perl
MongoDB driver~\cite{perlMongo}. With
\begin{lstlisting}
  polytope > db_info();
DATABASE: LatticePolytopes
This database contains various classes of lattice polytopes.

Collection: SmoothReflexive
A complete collection of smooth reflexive lattice polytopes in dimensions up to 9, up to lattice equivalence. [...]
\end{lstlisting}
we can query which collection groups are available. The
collection group and collection we want to use for our search are then
specified with the keywords \texttt{db} and \texttt{collection} in any
access function. The query itself is given as a perl hash. The query
is not processed by \polymake\ but directly handed over to MongoDB, so
it allows all queries specified in the MongoDB query language. A
specification of the full query language and its use from within perl
can be found in the documentation of MongoDB~\cite{mongodb} and the
perl driver for it~\cite{perlMongo}.
 
Here is an example returning an array of results.
\begin{lstlisting}
polytope > $parray=db_query({"DIM"=>3, "N_FACETS"=>5}, 
polytope(2) >                db=>"LatticePolytopes",
polytope(3) >                collection=>"SmoothReflexive");
polytope > print $parray->size;
4
\end{lstlisting}
This shows that there are $4$ polytopes in the collection
\textit{SmoothReflexive} that have dimension $3$ and $5$ facets. Using
a loop over this array or a database cursor we can check properties of
each object returned. For example
\begin{lstlisting}
polytope > $cursor=db_cursor({"DIM"=>3, "N_FACETS"=>5}, 
polytope(2) > db=>"LatticePolytopes",
polytope(3) > collection=>"SmoothReflexive");
polytope > while ( !$cursor->at_end() ) { 
polytope(2) > $p=$cursor->next(); 
polytope(3) > print $p->N_LATTICE_POINTS, " "; 
polytope(4) > }
34 30 31 30
\end{lstlisting}
defines a cursor over the collection \textit{SmoothReflexive}
successively returning all polytopes that satisfy the restrictions
given in the query, \textit{i.e.}, that have $5$ facets in dimension
$3$. Here it tells us that among the four polytopes found above, two
have $30$, one has $31$ and one has $34$ lattice points.

\section{Decomposing Smooth Fano Polytopes}

We illustrate the use of \polydb\ and its interface to \polymake\ with
a computation that uses the the collection \textit{SmoothReflexive} in
the collection group \textit{LatticePolytopes} to compute
decompositions of smooth Fano polytopes in dimensions $1$ to $8$. With
our computations we start a new statistics that counts how many of the
smooth Fano polytopes can be generated from lower dimensional smooth
Fano polytopes with some simple known polytope construction method
that preserves both smoothness and reflexiveness of the polytope. We
consider three methods in this paper and determine how many of the
smooth Fano polytopes in these dimensions are
\begin{itemize}
\item free sums of two smooth Fano polytopes
\item a smooth skew bipyramid over a smooth Fano polytope as defined
  in~\cite{AJP}, or
\item a generalized simplex sum of a smooth Fano polytope with a
  smooth simplex. This new construction method will be defined below.
\end{itemize}
All smooth skew bipyramids and many of the free sums are also
generalized smooth simplex sums. We will also provide the total number
of smooth Fano polytopes that can be decomposed with at least one of
these constructions. The results are collected in
Table~\ref{tab:sr_prod}.

We briefly explain the relevant notions. More background can,
\textit{e.g.}, be found in the book of Ewald~\cite{Ewald96}.  Let
$P\subseteq\R^d$ be a polytope with vertices
$v_1, \ldots, v_r\in \R^d$, \textit{i.e.},
\begin{equation}
  P\ := \ \conv(v_1, \ldots, v_r)\label{eq:interior}
\end{equation}
is the convex hull of these points and none of the $v_i$ can be
omitted in the definition. We assume that $P$ is full dimensional,
\textit{i.e.}, the affine hull of $P$ is $\R^d$ (otherwise we can pass
to a subspace). A polytope can equally be given as the intersection of
a finite number of half-spaces in the form
\begin{equation}
  P\ = \ \{\; x\mid Ax\; \le\; b\;\}\label{eq:exterior}
\end{equation}
for some $A\in\R^{s\times d}$ and $b\in\R^s$. We can again assume
that no inequality is redundant in this definition. In this case the
rows of $A$ are the \emph{facet normals} of $P$. A \emph{facet} $F$ of
$P$ is the set of all $x\in P$ that satisfy one of the inequalities
in \eqref{eq:exterior} with equality. A \emph{face} of $P$ is the
common intersection in $P$ of a subset of the facets (this may be
empty). The vertices, which are the faces of dimension $0$, are in the
common intersection of at least $d$ facets.

If $0$ is strictly contained in the interior of $P$, then the
\emph{polar} or \emph{dual polytope} is defined as
\begin{equation*}
  P^\vee\ :=\ \{\; v\mid \langle v, x\rangle\le 1\ \text{ for all }\ x\in P\;\}\,.
\end{equation*}
In fact, a finite subset of the inequalities in this definition
suffice to define $P^\vee$ (those corresponding to the vertices of
$P$), so that $P^\vee$ is again a polytope. Further, we have
$(P^\vee)^\vee=P$.

A lattice $\Lambda$ is the integral span of a linearly independent set
of vectors in $\R^n$. Up to a linear transformation we can assume that
$\Lambda$ is the integer lattice $\Z^d\subset\R^n$, and by passing to
a subspace we can assume that $n=d$. With these assumptions a polytope
$P$ is a \emph{lattice polytope} if all its vertices are in $\Z^d$.

In this case we can assume that both $A$ and $b$ are integral in
\eqref{eq:exterior}, and that the greatest common divisor of the
entries of each row of $A$ (\textit{i.e.}, of th entries of each facet
normal) is $1$.  A lattice polytope $P$ is \emph{reflexive} if
$P^\vee$ is again a lattice polytope. In this case $b=\1$ in
\eqref{eq:exterior}, and for both $P$ and $P^\vee$ the origin is the
unique interior lattice point. $P$ is \textit{smooth} if the vertices
of any facet of $P$ are a lattice basis of $\Z^d$. In this case $0$ is
a strictly interior point of $P$ and each facet has exactly $d$
vertices, so $P$ is \emph{simplicial}. Moreover, the polar polytope is
again a lattice polytope (in the dual lattice) whose vertices are the
facet normals (the rows of $A$), so $P$ is also reflexive.  Note that
in the literature sometimes the polytopes polar to the ones defined
here are called smooth.

\begin{figure}[t]
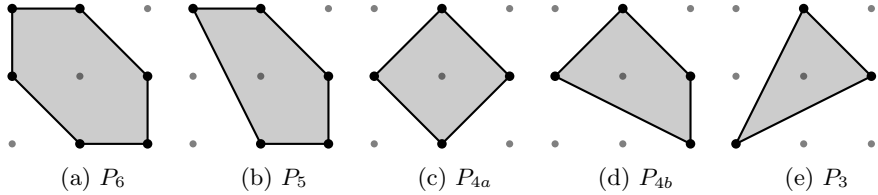

  \smoothFano
  \caption{The five $2$-dimensional smooth Fano polytopes.}
  \label{fig:2dmin_Fano}
\end{figure}

It follows from a result of Hensley~\cite{Hensley1983} and Lagarias
and Ziegler~\cite{MR1138580} that there are only finitely many smooth
reflexive polytopes in each dimension up to lattice equivalence
(affine transformations preserving $\Z^d$), as reflexive polytopes
have exactly one interior lattice point. See
Figure~\ref{fig:2dmin_Fano} for the list of such polytopes in
dimension $2$. The complete list is contained in \polydb\ for $d\le 9$
in the collection \textit{SmoothReflexive} of the database
\textit{LatticePolytopes}. Note however, that in the database we
follow the above mentioned alternative definition and list the duals
of the ones defined here. Yet, for the purpose of the following
constructions it is easier to work with the definition given above, so
we will use that one in the following. This requires that in the
scripts we use for our computations below we polarize the polytopes
obtained from the database. Sometimes this is, however, only done
implicitly. This is saves computation time, as it follows from the
design of \polymake\ that for reflexive polytopes the facets of the
polytope are the vertices of its dual. Also, as we will see below,
most constructions can also be given for the duals of the polytopes.

\begin{figure}[t]
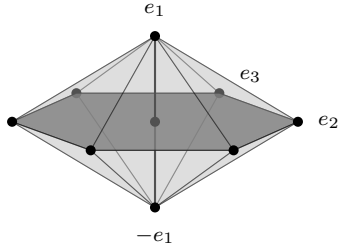
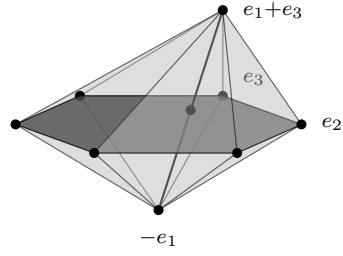
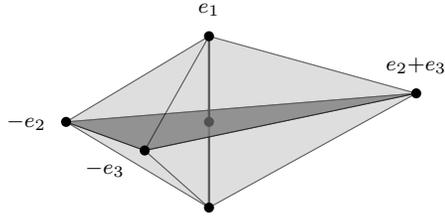
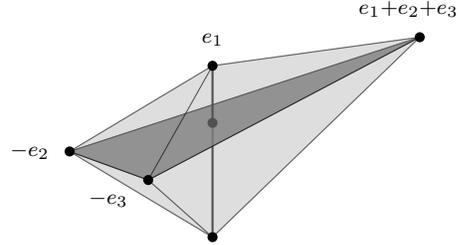

  \begin{subfigure}{.45\linewidth}
    \properBip 
    \caption{The free sum of a hexagon and a segment. This is at the
      same time also a proper bipyramid over the hexagon.}
    \label{subfig:properbip}
  \end{subfigure}
  \hfill
  \begin{subfigure}{.45\linewidth}
    \skewBip
    \caption{A skew bipyramid of a hexagon. The top apex has been shifted to $e_1+e_3$.\newline\mbox{ }}
    \label{subfig:skewbip}
  \end{subfigure}
  \begin{subfigure}{.495\linewidth}
    \skewSimplexSum
    \caption{A generalized simplex sum of a segment and a triangle.}
    \label{subfig:skewbip1}
  \end{subfigure}
  \hfill
  \begin{subfigure}{.495\linewidth}
    \skewSimplexSumTwo
    \caption{Another generalized simplex sum of a segment with a triangle.}
    \label{subfig:skewbip2}
  \end{subfigure}
  \caption{Polytope constructions}
  \label{fig:polytopeconstructions}
\end{figure}
We introduce several methods to construct a smooth Fano polytope from
smaller ones.  The most well known construction is the \emph{free sum}
of two polytopes $P\subseteq \R^a$ and $Q\subseteq\R^b$ that both
contain the origin in their interior. This is the polytope
\begin{equation*}
  P\oplus Q\ :=\ \conv\left(\left\{(v,0)\in\R^{a+b}\mid v\in P\right\}\cup\left\{(0,w)\in\R^{a+b}\mid w\in Q\right\}\right)\,.
\end{equation*}
We can also define this on the dual side. The \emph{product} of
polytopes $P$ and $Q$ is the polytope
\begin{equation*}
  P\times Q\ := \ \{\;(x,y)\mid x\in P,\; y\in Q\;\}\;.
\end{equation*}
Then, if $P$ and $Q$ contain the origin in their interior,
\begin{equation}
  P\oplus Q\,=\,(P^\vee \times Q^\vee)^\vee\, .\label{eq:dualsum}
\end{equation} 
We will use this dual
definition for the detection of free sums among the smooth Fano
polytopes. See Figure~\ref{subfig:properbip} for an example.

A \emph{bipyramid} over a polytope $P$ is the free sum of $P$ with a
segment $S$ containing the origin in the interior. More generally, we
say that $Q$ is a \emph{skew bipyramid} over $P$ if $Q$ has the same
combinatorial type (the same face lattice) as a bipyramid over
$P$. The two vertices coming from vertices of $S$ are the two
\emph{apices} of $Q$.

If $P$ is a smooth $d$-dimensional Fano polytope then we call the free
sum with the segment $[-1,1]$ the \emph{smooth bipyramid} over
$P$. Let $v$ be a vertex of $P$ and $\bar v$ its embedding into
$\R^{d+1}$ by adding a $0$ at the end. Then the \emph{smooth skew
  bipyramid} for vertex $v$ as defined by Assarf et al.~\cite{AJP} is
the polytope
\begin{equation*}
  \sbip(P,v)\ :=\ \conv\left(\, P\times \{0\} \cup \left\{-e_{d+1}, \bar v+e_{d+1}\right\}\right)\,.
\end{equation*}
Figure~\ref{subfig:skewbip} shows an example of this definition. More
generally, we say that $Q$ is a \emph{smooth generalized skew
  bipyramid} over $P$ if $Q$ is a skew bipyramid over $P$ such that
the two apices have lattice distance $1$ from $P$. This class contains
all smooth bipyramids and smooth skew bipyramids. The following
proposition is an extension of Lemmas $1$, $2$ and $3$ of
\cite{AJP}. The proof easily carries over into this more general
setting.
\begin{proposition}
  Let $P$ and $Q$ be smooth Fano polytopes. Then the free sum
  $P\oplus Q$, the smooth bipyramid and any smooth (generalized) skew
  bipyramid over $P$ are again smooth Fano polytopes.\hfill\mbox{ }\qed
\end{proposition}
We further generalize this construction. Let $P\subseteq\R^a$ be a
smooth Fano polytope and $Q\subseteq \R^b$ a smooth Fano simplex (this
is unique up to lattice equivalence). Let $v$ be a vertex of $Q$. Then
$R:=P\oplus Q$ is a smooth Fano polytope and also any polytope $R'$
obtained from $R$ by replacing $v$ with a lattice point $v'$ in the
hyperplane $\R^a+v\subseteq \R^{a+b}$, as long as $R$ and $R'$ have the
same combinatorial type. This is again a simple extension of the
proposition above. We call those polytopes \emph{smooth generalized
  simplex sums}. Figures~\ref{subfig:skewbip1}
and~\ref{subfig:skewbip2} show two examples. Observe that any smooth
(generalized skew) bipyramid is a simplex sum.

We can use \polymake\ and \polydb\ to detect all free sums and smooth
generalized simplex sums among the smooth Fano polytopes.  Clearly,
these two constructions overlap in various ways. Any proper bipyramid
over all polytope $P$ is also the free sum of a $P$ with a segment, and
polytopes may have more than one possible decomposition into a free
sum. Many of the various possibilities to place the vertex $v'$ for a
generalized smooth simplex sum are lattice equivalent. Our approach to
detect all different instances is as follows: For a fixed dimension
$d$ we consider all possible splits of $d$ as a sum of dimensions $a$
and $b$ and compute all free sums of smooth Fano polytopes in these
two dimensions and all simplex sums of an $a$-dimensional smooth Fano
polytope with a $b$-dimensional simplex. For each such polytope we run
through the list of $d$-dimensional smooth Fano polytopes, check for
lattice equivalence and store the name of the polytope we have
found. We could also store the way we obtained it alongside, so that
in the end we have a list of all possible splits for a given
$d$-dimensional smooth Fano polytope.
\begin{figure}[t]
  \usebox{\productscript}
  \caption{A function to detect all free sums among the smooth Fano
    polytopes. This is a shortened version of the script given at
    \cite{polymakeSmoothFano}.}
  \label{fig:ident_prod}
\end{figure}

\begin{table}[t]
  \caption{Free sums, skew bipyramids and generalized smooth simplex sums among the smooth Fano polytopes. The rows denoted by \textit{sg simplex-n sums} for $n$ between $1$ and $8$ count the simplex sums with a simplex of dimension $n$. The row denoted by \textit{total sg simplex sums} gives the number of different generalized smooth simplex sums with a simplex of any dimension. The row \textit{total decomposable} counts the number of different polytopes among the free sums and the smooth generalized simplex sums.}
  \label{tab:sr_prod}
  \centering
  \begin{tabular}{rrrrrrrr}
    \toprule
    dimension&2&3&4&5&6&7&8\\
    \midrule
    smooth Fano polytopes&5&18&124&866&7622&72256&749892\\
    \midrule
    free sums&1&5&28&176&1361&11760&112285\\
    skew bipyramids&1&9&57&489&4323&43777&466770\\
    \midrule
    sg simplex-1 sums&2&13&66&556&4700&47076&495092\\
    sg simplex-2 sums&1&3&31&232&2403&25157&284249\\
    sg simplex-3 sums&-&1&4&52&515&6635&83730\\
    sg simplex-4 sums&-&-&1&5&81&961&14598\\
    sg simplex-5 sums&-&-&-&1&6&114&1609\\
    sg simplex-6 sums&-&-&-&-&1&7&155\\
    sg simplex-7 sums&-&-&-&-&-&1&8\\
    sg simplex-8 sums&-&-&-&-&-&-&1\\
    \midrule
    total sg simplex sums&3&16&93&708&6283&61961&657380\\
    \midrule
    total decomposable &3&16&96&712&6346&62331&660792\\
    \bottomrule
  \end{tabular}
\end{table}
We did the computation up to dimension $8$. The results are given in
Table~\ref{tab:sr_prod}. The free sums can be obtained with the small
scripts given in Figure~\ref{fig:ident_prod}. The first script
\texttt{identify\_smooth\_polytope} takes a smooth Fano polytope,
identifies it in the database and returns its name. The identification
is based on the \polymake\ function
\texttt{lattice\_isomorphic\_smooth\_polytopes}, that reduces the
check whether two lattice polytopes are lattice isomorphic to a
colored graph isomorphism problem (which is solved using
\texttt{bliss}~\cite{bliss} or \texttt{nauty}~\cite{nauty}). Note that
there is also the extension
\textit{LatticeNormalization}~\cite{latticeNormalization} to
\polymake\ that computes the lattice normal form of a lattice polytope
(see~\cite{1301.6641} for a definition), but the reduction to colored
graph isomorphism is more efficient for smooth polytopes. The simpler
problem of checking combinatorial isomorphisms (\textit{i.e.}, graph
isomorphism) can also be done with the \polymake-function
\texttt{canonical\_hash} (also based on \texttt{bliss} or
\texttt{nauty}). The second function \texttt{all\_free\_sums\_in\_dim}
computes all possible free sums that lead to a $d$-dimensional smooth
Fano polytope. As the database contains the polytopes dual to the ones
we consider we use~\eqref{eq:dualsum} and compute products instead of
sums to avoid explicit dualization. For each product the function
calls \texttt{identify\_smooth\_polytope} to identify it in the
database. The function returns a list of all names (\texttt{\_id}s)
found in this way. If \texttt{splitinfo} is set to $1$ it also returns
all pairs of summands.

For the computation of the smooth generalized simplex sums we used the
function \texttt{all\_skew\_simplex\_sums\_in\_dim} available
at~\cite{polymakeSmoothFano}. For each combination of an
$a$-dimensional smooth Fano polytope and a $b$-dimensional simplex
with $d=a+b$ we compute all possible lattice points for the shifted
vertex $v'$, construct the polytope and again use
\texttt{identify\_smooth\_polytope} to identify it in the
database. Computation of all possible $v'$ requires the computation of
all lattice points in the hyperplane $\R^a+v\subseteq R^{a+b}$ that
lead to a lattice polytope with the same combinatorial type as the
proper free sum. This can be reduced to enumerating lattice points in
the interior of a polytope, which is done in \polymake\ via the
interface to \texttt{Normaliz}~\cite{Normaliz312}. As above the
function returns a list of \texttt{id}s, and also all possible
decompositions into a simplex sum if \texttt{splitinfo} is set to $1$.
\begin{table}[t]
  \begin{subtable}{.34\linewidth}
    \small
      \begin{tabular}{>{\!}r>{\!}r>{\!}r>{\!}r>{\!}r}
             $0$&  $0$& $ 0$& $ 0$& $1$\\
             $0$&  $0$& $ 1$& $ 0$& $-1$\\
             $0$&  $0$& $-1$& $ 0$& $0$\\
             $0$&  $0$&  $0$& $-1$& $0$\\
             $0$&  $0$&  $0$&  $0$& $-1$\\
             $-1$& $ 0$& $ 0$& $ 0$&$ 0$\\
             $0$& $-1$&  $0$&  $0$& $0$\\
             $0$&  $1$&  $0$&  $0$& $1$\\
             $1$&  $0$&  $0$&  $1$& $2$
      \end{tabular}
    \medskip
    \caption{Dual of \texttt{F.5D.0116}\newline\mbox{ }\newline\mbox{ }}
    \label{tab:f0116}
  \end{subtable}
  \hspace*{3cm}
  \begin{subtable}{.34\linewidth}
    \small
    \begin{tabular}{>{\!}r>{\!}r>{\!}r>{\!}r>{\!}r}
      \cellcolor{lightgray}$ 0$&\cellcolor{lightgray} $ 1$ &\cellcolor{lightgray} $ 1$&$0$&$0$\\
      \cellcolor{lightgray}$-1$&\cellcolor{lightgray} $ 0$ &\cellcolor{lightgray} $ 0$&$0$&$0$\\
      \cellcolor{lightgray}$ 0$&\cellcolor{lightgray} $ 0$ &\cellcolor{lightgray} $ 1$&$0$&$0$\\
      \cellcolor{lightgray}$ 0$&\cellcolor{lightgray} $-1$ &\cellcolor{lightgray} $ 0$&$0$&$0$\\
      \cellcolor{lightgray}$ 0$&\cellcolor{lightgray} $ 0$ &\cellcolor{lightgray} $-1$&$0$&$0$\\
      \cellcolor{lightgray}$ 1$&\cellcolor{lightgray} $ 0$ &\cellcolor{lightgray} $-1$&$0$&$0$\\
      $ 0$& $ 0$ & $ 0$&$-1$&$ 0$\\
      $ 0$& $ 0$ & $ 0$&$ 0$&$-1$\\
      $ 0$& $ 0$ & $-2$&$ 1$&$ 1$\\
    \end{tabular}

    \medskip
    \caption{Dual of \texttt{F.3D.0112} extended with skew triangle}
    \label{tab:f0112}
  \end{subtable}
 
  \begin{subtable}{.34\linewidth}
    \small
    \begin{tabularx}{.2\linewidth}{>{\!}r>{\!}r>{\!}r>{\!}r>{\!}r}
      \cellcolor{lightgray}$-1$&\cellcolor{lightgray} $ 0$ &\cellcolor{lightgray} $ 0$ &\cellcolor{lightgray} $ 0$&$0$\\
      \cellcolor{lightgray}$ 0$&\cellcolor{lightgray} $ 0$ &\cellcolor{lightgray} $ 0$ &\cellcolor{lightgray} $-1$&$0$\\
      \cellcolor{lightgray}$ 0$&\cellcolor{lightgray} $ 1$ &\cellcolor{lightgray} $ 0$ &\cellcolor{lightgray} $ 1$&$0$\\
      \cellcolor{lightgray}$ 0$&\cellcolor{lightgray} $-1$ &\cellcolor{lightgray} $ 0$ &\cellcolor{lightgray} $ 0$&$0$\\
      \cellcolor{lightgray}$ 0$&\cellcolor{lightgray} $ 0$ &\cellcolor{lightgray} $ 0$ &\cellcolor{lightgray} $ 1$&$0$\\
      \cellcolor{lightgray}$ 0$&\cellcolor{lightgray} $ 0$ &\cellcolor{lightgray} $-1$ &\cellcolor{lightgray} $ 0$&$0$\\
      \cellcolor{lightgray}$ 1$&\cellcolor{lightgray} $ 0$ &\cellcolor{lightgray} $ 1$ &\cellcolor{lightgray} $ 2$&$0$\\
      $ 0$& $ 0$ & $ 0$&$ 0$&$-1$\\
      $ 0$& $ 0$ & $ 0$&$-1$&$ 1$\\
    \end{tabularx}
    \medskip
    \caption{Dual of \texttt{F.4D.0008} extended with segment}
    \label{tab:f0008}
  \end{subtable}
  \hspace*{3cm}
  \begin{subtable}{.34\linewidth}
    \small
    \begin{tabularx}{.2\linewidth}{>{\!}r>{\!}r>{\!}r>{\!}r>{\!}r}
      \cellcolor{lightgray}$-1$&\cellcolor{lightgray} $ 0$ &\cellcolor{lightgray} $ 0$ &\cellcolor{lightgray}$ 0$&$0$\\
      \cellcolor{lightgray}$ 0$&\cellcolor{lightgray} $ 0$ &\cellcolor{lightgray} $ 0$ &\cellcolor{lightgray}$-1$&$0$\\
      \cellcolor{lightgray}$ 0$&\cellcolor{lightgray} $ 1$ &\cellcolor{lightgray} $ 0$ &\cellcolor{lightgray}$-1$&$0$\\
      \cellcolor{lightgray}$ 0$&\cellcolor{lightgray} $-1$ &\cellcolor{lightgray} $ 0$ &\cellcolor{lightgray}$ 0$&$0$\\
      \cellcolor{lightgray}$ 0$&\cellcolor{lightgray} $ 0$ &\cellcolor{lightgray} $ 0$ &\cellcolor{lightgray}$ 1$&$0$\\
      \cellcolor{lightgray}$ 0$&\cellcolor{lightgray} $ 0$ &\cellcolor{lightgray} $-1$ &\cellcolor{lightgray}$ 0$&$0$\\
      \cellcolor{lightgray}$ 1$&\cellcolor{lightgray} $ 0$ &\cellcolor{lightgray} $ 1$ &\cellcolor{lightgray}$ 2$&$0$\\
      $ 0$& $ 0$ & $ 0$&$ 0$&$-1$\\
      $ 0$& $ 0$ & $ 0$&$ 1$&$ 1$\\
    \end{tabularx}
    \medskip
    \caption{Dual of \texttt{F.4D.0019} extended with segment}
    \label{tab:f0019}
  \end{subtable}

  \caption{Simplex sums leading to the dual of \texttt{F.5D.0116}}
  \label{tab:decomp}
\end{table}
You can save the scripts to a file in the current folder and load this
into \polymake\ via
\begin{lstlisting}
  polytope> script(<filename>);
\end{lstlisting}
Then the classification, \textit{e.g.} in dimension 4, is obtained
with
\begin{lstlisting}
  polytope > $fs = all_free_sums_in_dim(4);
  polytope > print $fs->size;
  28
  polytope > $sb = all_skew_bipyramids_in_dim(4);
  polytope > print $fs->size;
  57
  polytope > $s = new Set<String>;
  polytope > foreach (1..4) {
  polytope(2) > $st = skew_simplex_sums_in_dim(4,$_);
  polytope(3) > print $st->size, " ";
  polytope(4) > $s += $st;
  polytope(5) > }
  66 31 4 1
  polytope > print $s->size;
  93
  polytope > print (($fs+$s)->size);
  96
\end{lstlisting}
The scripts containing the functions
\texttt{all\_skew\_bipyramids\_in\_dim} for skew bipyramids and
\texttt{skew\_simplex\_sums\_in\_dim} for generalized smooth simplex
sums are available from~\cite{polymakeSmoothFano} and allow to store
the possible decompositions. Note that the computation time for the
decompositions grows quickly in the dimension. While dimension $4$
runs in a few minutes on an Intel Xeon E5-4650, computations in
dimension $8$ took over a month.

From these computations we can, \emph{e.g.}, see that we have three
different decompositions of the dual of the $5$-dimensional polytope
with index \texttt{F.5.0116}. Its vertices are the rows of the matrix
in Table~\ref{tab:f0116}. We can decompose this into three different
simplex sums. One is over dual of the $3$-dimensional polytope $P_3$
with index \texttt{F.3D.0112}. This is shown in Table~\ref{tab:f0112},
where the shaded part corresponds to the vertices of $P_3$. The
shifted vertex of the triangle is $[0,0,-2,1,1]$. Note that the
vertices are given as obtained by dualization from the
database. Hence, the equality is not directly visible from the
vertices, as the two polytopes differ by a lattice isomorphism. We can
check this with polymake.
\begin{lstlisting}
  polytope > $p3_ext = new Polytope(VERTICES=>
polytope(2) > [[1,0,1,1,0,0],[1,-1,0,0,0,0],[1,0,0,1,0,0],
polytope(3) > [1,0,-1,0,0,0],[1,0,0,-1,0,0],[1,1,0,-1,0,0],
polytope(4) > [1,0,0,0,-1,0],[1,0,0,0,0,-1],
polytope(5) > [1,0,0,-2,1,1]]);
polytope > $p5 = new Polytope(VERTICES=>
polytope(2) > [[1,0,0,0,0,1],[1,0,0,1,0,-1],[1,0,0,-1,0,0],
polytope(3) > [1,0,0,0,-1,0],[1,0,0,0,0,-1],[1,-1,0,0,0,0],
polytope(4) > [1,0,-1,0,0,0],[1,0,1,0,0,1],[1,1,0,0,1,2]]);
polytope > print lattice_isomorphic_smooth_polytope(
polytope(2) > polarize($p3_ext),polarize($p5));
1
\end{lstlisting}
Here, the variable \texttt{\$p3\_ext} contains the polytope $P_3$ and
\texttt{\$p5} is $P_5$. As above we need to dualize for the
isomorphism check. The check returns $1$, which is the
\texttt{true}-value for \texttt{polymake}.

The other two decompositions are over the $4$-dimensional polytopes
$P_4^1$ and $P_4^2$ with index \texttt{F.4D.0008} and
\texttt{F.4D.0019}. Those are shown in Tables~\ref{tab:f0008}
and~\ref{tab:f0019}. Again, the vertices of $P_4^1$ and $P_4^2$ are
shaded. The shifted vertices of the $1$-dimensional simplex are in the
last line.

With this simple computation we have seen that over $80\%$ of the
smooth Fano polytopes can be obtained from at least one of the
constructions considered here. Hence, for a structural description of
all smooth Fano polytopes it suffices to look at the remaining less
than $20\%$.

\end{document}